\documentclass{Tran-l}

 \newtheorem{thm}{Theorem}[section]
 \newtheorem{cor}[thm]{Corollary}
 \newtheorem{lem}[thm]{Lemma}
 \newtheorem{prop}[thm]{Proposition}
 \theoremstyle{definition}
 \newtheorem{defn}[thm]{Definition}
 \theoremstyle{remark}
 \newtheorem{rem}[thm]{Remark}
 
 \newcommand{\CC}{\mathbb{C}}

 \newcommand{\VS}{\mbox{VSH}}
 \newcommand{\YY}{\mathbb{Y}}
 
 \newcommand{\XX}{\mbox{$\mathbb{X}$}}

 \newcommand{\Sec}{\mathrm{Sec}}

  \newcommand{\PP}{\mbox{$\mathbb{P}$}}

  \def\sqw{\hbox{\rlap{\leavevmode\raise.3ex\hbox{$\sqcap$}}$%
\sqcup$}}

\def\cqfd{\ifmmode\sqw\else{\ifhmode\unskip\fi\nobreak\hfil
\penalty50\hskip1em\null\nobreak\hfil\sqw
\parfillskip=0pt\finalhyphendemerits=0\endgraf}\fi}

  \newtheorem{Teo}{Theorem}[section]

  \font\Dimfont=cmr8

  \font\Dimfont=cmr8

  \newtheorem{ProEs}[Teo]{Example}

  \newtheorem{ProRemark}[Teo]{Remark}

  \newcommand{\sexp}{s_\mathrm{exp}}

  \newcommand{\smin}{s_\mathrm{min}}

  \renewcommand{\P}{{\mathbb P}} 
\usepackage{xypic}

\begin{document}

\title{Codimension one decompositions and Chow varieties}

\author{E. Carlini}

\address{Department of Mathematics, University of Pavia, Via Ferrata 1, 27100 Italy}

\email{carlini@dimat.unipv.it}

\thanks{The author is a member of the PRIN01 project ``Spazi di moduli e teoria di Lie'' of the University of Pavia.}


\keywords{Waring's problem, polynomials decompositions, secant varieties, Chow varieties}




\begin{abstract}
A presentation of a degree $d$ form in $n+1$ variables as the sum of homogenous elements ``essentially'' involving $n$ variables is called a {\em codimension one decomposition}. Codimension one decompositions are introduced and the related Waring Problem is stated and solved. Natural schemes describing the codimension one decompositions of a generic form are defined. Dimension and degree formulae for these schemes are derived when the number of summands is the minimal one; in the zero dimensional case the scheme is showed to be reduced. These results are obtained by studying the Chow variety $\Delta_{n,s}$ of zero dimensional degree $s$ cycles in $\PP^n$. In particular, an explicit formula for $\deg\Delta_{n,s}$ is determined. 
\end{abstract}

\maketitle
\section{Introduction}

Usually, a homogeneous element in the polynomial ring $S=\CC[X_0,\ldots,X_n]$ is presented as a sum of monomials. In other words, we use the homogeneous structure to choose a vector space basis in each homogenous piece $S_d$ of $S$. Actually, we may want to write down $f\in S_d$ in different ways and this can be done also without selecting a vector space basis. Sum of powers decompositions (see \cite{IaKa}) are just an example:
\[f=l_1^d+\ldots l_s^d, l_1,\ldots,l_s\in S_1.\]

Sum of powers presentations have been widely studied classically in an attempt to produce a classification of homogeneous polynomials. The idea was to mimic what happens for quadratic forms and diagonalization (this has not been very successful, e.g. no effective algorithm is known to perform a sum of powers decomposition if $d>3$ and $n>1$). Nevertheless, information can be obtained on $f$ and on its zero locus by studying properties of its sum of powers decompositions (see \cite{Sal1}, e.g., page 252). In particular, it is useful to know how many summands are needed.

For a generic element $f\in S_d$, a parameters count shows that at least $\lceil{1 \over n+1 }{n+d\choose n}\rceil$ summands are needed for a sum of powers presentation of $f$ and before Clebsch's paper this number was believed to be always enough. In \cite{Clb}, it is shown that a generic ternary quartic (i.e. $n=2$, $d=4$) is the sum of $6$ and not of $5$ powers of linear forms. The existence of defective cases makes Waring Problem for forms so interesting:
\begin{quote}
{\em For each pair $(n,d)$ determine the minimal number of summands appearing in the sum of powers decompositions of the generic form of degree $d$ in $n+1$ variables.}
\end{quote}
{\em Defective pairs}, such as Clebsch's, were readily discovered, but the problem remained unsolved for almost one century. The complete answer was only recently found by Alexander and Hirschowitz (\cite{AH95}):
\begin{thm}\label{AHthm} A generic form of degree $d$ in $n+1$ variables is the sum of $s=\lceil{1 \over n+1}{n+d \choose d}\rceil$ powers of linear forms, unless
\begin{itemize}
\item $d=2$, where $s=n+1$ instead of $\lceil{n+2\over 2}\rceil$, or
\item $d=4$ and $n=2,3,4$, where $s=6,10,15$ instead of $5,9,14$ respectively, or
\item $d=3$ and $n=4$, where $s=8$ instead of $7$.
\end{itemize}
\end{thm}

As $l^d$ is a homogeneous element in the univariate ring $\CC[l]$, a sum of powers decomposition can be viewed as a presentation of a form as sum of forms ``essentially'' involving one variable. With this in mind is natural to consider other presentations of this kind, e.g. binary decompositions, where the summands essentially involve two variables (see \cite{Ca04} and \cite{Ca04JA}). We can also move to the opposite end of the spectrum and consider {\em codimension one decompositions}, where the summands essentially involve one variable less than the original form.

The study of codimension one decompositions is the object of this paper. In section \ref{waring}, we will address and solve the analogous Waring type problem obtaining the following results:

\medskip\noindent
{\bf Theorem \ref{sminTHM}. }{\em The generic form of degree $d$ in $n+1$ variables is the sum of $\min\{s: ns-{d-s+n\choose n }\geq 0\}$ codimension one forms and no fewer.}

\medskip\noindent
{\bf Corollary \ref{AHsimeq}. }{\em Let $n\geq 2$. The minimal number of summands appearing in the codimension one decompositions of the generic form of degree $d$ in $n+1$ variables is the expected one (see Definitions \ref{defsmin} and \ref{defsexp} ) if and only if
\[d=2,3\mbox{ for any } n\geq 2\]
or
\[d=4,5,6\mbox{ and }8 \mbox{ for }n=2.\]}
\medskip

In particular, the Corollary shows how codimension one and sum of powers decompositions are deeply different. In the sum of powers case the expected number of summands almost always works with some exceptions, as shown in Theorem \ref{AHthm}. But, for codimension one decompositions, exactly the opposite happens: only in some cases the expected number of summands works and almost all the pairs are defective.

In section \ref{VSH}, we introduce and study a natural scheme, $\VS$ (see Definition \ref{defVS}), describing the codimension one decompositions of a generic form. This is mostly done in the spirit of \cite{RS00} and we obtain dimension and  degree formulae:

\medskip\noindent
{\bf Theorem \ref{DIMdegTHM}. }{\em Let $F\subset\PP^n$ be a generic degree $d$ hypersurface and let $s=\smin(n,d)$, then
\begin{itemize}
\item $\dim\VS(F,s)=ns-{d-s+n\choose n}$;
\item $\deg \VS(F,s)={ns-1\choose n-1}\cdot{n(s-1)-1\choose n-1}\cdot\ldots\cdot 1$.
\end{itemize}}
\medskip

In the zero dimensional cases we also get reducedness:

\medskip\noindent
{\bf Proposition \ref{PropRED}. }{\em Let $F\subset\PP^n$ be a generic degree $d$ hypersurface and let $s=\smin(n,d)$. If $\VS(F,s)$ is zero dimensional, then it is reduced.}
\medskip

These results are obtained through a careful study of the Chow variety $\Delta_{n,s}$ and of some special linear sections of this. In particular, we obtain a description of the tangent space in a generic point and an explicit degree formula (see Proposition \ref{degDIMprop}).

I wish to thank Aldo Conca and Anthony Geramita for their help with the algebraic claim in the proof of Proposition \ref{degDIMprop}: the latter for giving me an idea of a proof and the former for showing me a much better proof than the one I had. The referee's comments were extremely useful. In particular, I wish to thank her/him for showing to me the connection with Chow varieties. The hospitality of the Mathematics Department of the University of Genoa and the financial support of the Mathematics Department of the University of Pavia were appreciated.

{\bf Notation:} we work with the polynomial ring $S=\CC[X_0,\dots,X_n]$ and its ring of differential operators $T=\CC[\partial_0,\dots,\partial_n]$, i.e. $\partial_i$ acts as the partial derivation ${\partial\over\partial x_i}$. In particular, $S_1$ and $T_1$ are dual to each other and we let $\PP^n=\PP T_1$ and $\check\PP^n=\PP S_1$. A form $f\in S_d$ defines a hypersurface $F=V(f)\subset\PP^n$ and linear spaces $(f^\perp)_s\subset\PP T_s$ which we will denote as $F^\perp\subset\PP T_s$ with abuse of notation.
We work over the complex number field $\CC$, but any algebraically closed field of characteristic 0 could be used instead (in positive characteristic problems arise because of the coefficients produced by differentiating).

\section{Waring Problem for codimension one decompositions}\label{waring}

In what follows we need some basic facts about {\em apolarity theory} (see \cite{Ge}). In particular, we consider the polynomial rings $S=\CC[X_0,\dots,X_n]$ and $T=\CC[\partial_0,\dots,\partial_n]$ where $S$ has a $T$-module structure given by the differentiation action, which we denote with ``$\circ$''. Given a form $f\in S_d$, $f^\perp=\{D\in T : D\circ f=0\}\subset T$ denotes the homogeneous ideal of derivations annihilating $f$ and $T/f^\perp$ is an artinian Gorenstein ring with socle in degree $d$. Given a derivation $D\in T_d$, $D^{-1}=\{f \in S : D\circ f=0\}$ is a graduated sub-$T$-module. The apolarity pairing $S_d\times T_d\rightarrow\CC$ is perfect.
The link between apolarity and the study of polynomial decompositions is given by the classical Apolarity Lemma (for a proof see \cite{RS00}):
\begin{lem}[Apolarity Lemma]\label{ApolarityLemma} Let $f\in S_d$, then the following are equivalent:
\begin{enumerate}
\item $f=l_1^d+\dots + l_s^d$, where the $l_i$'s are pairwise non-proportional linear forms;
\item $f^\perp\supset I_{\XX}$, where $I_{\XX}$ is the ideal of the set of $s$ points $\XX=\{l_1,\dots,l_s\}\subset\check\PP^n$.
\end{enumerate}
\end{lem}

\begin{defn} Let $g\in S_d$, then $g$ is called a {\em codimension one form} if $(g^\perp)_1\neq 0$. Given a form $f\in S_d$, a {\em codimension one decomposition} of $f$ is a presentation
\[f=\hat f_1+\ldots \hat f_s\]
where $\hat f_i, i=1,\dots,s$, are codimension one forms of degree $d$.
\end{defn}
\begin{rem}
If $g\in S_d$ is a codimension one form, then there exists a linear change of variables $X_i\mapsto Y_i,i=0,\dots,n$, such that $g(Y)$ only involves $n$ variables, i.e. $g(Y)\in\CC[Y_1,\dots,Y_n]$.
\end{rem}
\begin{rem}
Codimension one forms can be nicely described in geometric terms. Let $\nu_d:\PP S_1=\check\PP^n\rightarrow \PP S_d$ be the $d$-uple embedding and denote by $\langle\nu_d(H)\rangle\subset\PP S_d$ the linear span of the image of a hyperplane $H\subset{\check\PP^n}$. Then the variety
\[\hat V_{n,d}=\bigcup_{H\in{\PP^n}}\langle\nu_d(H)\rangle\]
parameterizes the codimension one forms in $\PP S_d$.
Notice that $\hat V_{n,d}$ is a determinantal variety defined by the maximal minors of a $(n+1)\times{n+d-1\choose n}$ matrix of linear forms (which is also a catalecticant matrix, see \cite{Ge1}). But, in general, it is not standard determinantal and $\dim\hat V_{n,d}=n+{d+n-1\choose n-1}$.
\end{rem}

No algorithm is known to determine a codimension one decomposition of a given form, thus it is interesting to study quantitative aspects of such a presentation, e.g. the number of summands $s$. As in the case of sum of powers decompositions, we begin by studying the behavior of {\em generic} forms (see Remark \ref{generic}):
\begin{defn}\label{defsmin} $\smin(n,d)$ is the minimal number of summands appearing in the codimension one decompositions of the generic form of degree $d$ in $n+1$ variables.
\end{defn}
\begin{rem}\label{generic} To make clearer what we mean by {\em generic}, it is useful to recall some geometry. Let $\Sec_t(\hat V_{n,d})$ be the variety of $t+1$ secants, i.e. the {\em closure} of the union of the $\P^t$'s spanned by points of $\hat V_{n,d}$. Then, the generic $f\in\Sec_t(\hat V_{n,d})$ is a sum of $t+1$ codimension one forms. In these terms, $t=\smin(n,d)-1$ is the smallest integer such that $\Sec_t(\hat V_{n,d})=\PP S_d$. Studying the decomposition of {\em any} form one completely loses this nice geometric interpretation and the problem gets considerably harder (even in the sums of power case no complete solution is known!). In this paper we will only deal with the generic case.
\end{rem}
An estimate for $\smin$ can be easily determined. As
\[\dim \Sec_{s-1}(\hat V_{n,d})\leq s\dim \hat V_{n,d}+s-1\]
the condition $\Sec_{s-1}(\hat V_{n,d})=\PP S_d$ gives an inequality and solving it we get
\[\smin(n,d)\geq \left\lceil{1\over n+{d+n-1\choose n-1}}{d+n\choose n}\right\rceil.\]
\begin{defn}\label{defsexp} The {\em expected value} of $\smin(n,d)$ is $\sexp(n,d)=\lceil{1\over n+{d+n-1\choose n-1}}{d+n\choose n}\rceil$. If $\smin(n,d)\neq\sexp(n,d)$ the pair $(n,d)$ is said to be {\em defective}.
\end{defn}
The Waring Problem for codimension one decompositions can be stated as follows:
\begin{quote}
{\em For each pair $(n,d)$ determine the minimal number of codimension one forms needed for the decomposition of the generic form of degree $d$ in $n+1$ variables, i.e. compute $\smin(n,d)$.}
\end{quote}

Apolarity provides us with a strong tool to study codimension one decompositions:
\begin{lem}[Codimension One Lemma]\label{Cod1ApolarityLemma} Let $f\in S_d$, then the following are equivalent:
\begin{enumerate}
\item $f=\hat f_1+\dots + \hat f_s$, where the $\hat f_i$'s are codimension one forms of degree $d$ such that $(\hat f_i)^\perp_1\neq(\hat f_j)^\perp_1$ for $i\neq j$;
\item there exists $L_1\cdot\ldots\cdot L_s\in f^\perp$, where the $L_i$'s are pairwise non-proportional linear forms.
\end{enumerate}
\end{lem}
\begin{proof} If $f$ admits a codimension one decomposition with the property above, choose non-proportional linear forms $L_i\in(\hat f_i)^\perp_1, i=1,\dots,s$. Then $(L_1\cdot\ldots\cdot L_s)\circ f=0$ and the claim follows. Conversely, assume that there exists $L_1\cdot\ldots\cdot L_s\in f^\perp$ such that the hyperplanes $\{L_i=0\}, i=1,\dots,s$, are distinct and choose $N_{n,d}={d+n-1\choose n-1}$ generic points on each of them. If we denote by $\XX$ the resulting set of $sN_{n,d}$ points, then its defining ideal $I_{\XX}$ satisfies
\[(I_{\XX})_t\subseteq (f^\perp)_t, \mbox{ for } t>d.\]
But the inclusion also holds for $t\leq d$ by Bezout. Thus, the Apolarity Lemma yields
\[f=\sum_{i=1}^s\sum_{j=1}^{N_{n,d}}(l_{ij})^d\]
where $L_i\circ (\sum_{j=1}^{N_{n,d}}(l_{ij})^d)=0,i=1,\ldots,s,$ by construction. Hence the claim.
\end{proof}
\begin{rem}\label{sstar}
This Lemma produces a useful bound for $\smin$ almost without effort. Given a form $f\in S_d$ the ring $T/f^\perp$ is known to be artinian Gorenstein and to have socle in degree $d$. In particular, for a generic form, $f^\perp$ is generated in degree no smaller than $\lceil{d+1\over 2}\rceil$ and hence the Codimension One Lemma implies that $\smin(n,d)\geq\lceil{d+1\over 2}\rceil$.
\end{rem}
We will use Lemma \ref{Cod1ApolarityLemma} and basic Algebraic Geometry techniques to give an answer to the Waring Problem for codimension one decompositions.
\begin{defn}
Let $\Delta_{n,s}\subset\PP T_s$ be the variety of {\em totally decomposable forms} of degree $s$ in $n+1$ variables, i.e. a point of $\Delta_{n,s}$ represents a form which can be written as the product of $s$ linear forms.
\end{defn}
\begin{rem}\label{chowREM}
$\Delta_{n,s}$ is the Chow variety of zero-dimensional degree $s$ cycles in $\P^n$ and in these terms it has been widely studied in \cite{GKZ}. In particular, it is shown there how to find equations for $\Delta_{n,s}$ set-theoretically, but it is known that these equations do not generate the defining ideal.
\end{rem}
\begin{rem}
If we consider the symmetrization of the Segre product $\overbrace{\PP T_1\times\ldots\times\PP T_1 }^s$ and we embed it in $\PP T_s$, then we get $\Delta_{n,s}$. In particular, this shows that $\dim \Delta_{n,s}=ns$.
\end{rem}
With this setting we can develop a strategy to study codimension one decompositions in general. Given $f\in S_d$, we consider the smooth points of $\Delta_{n,s}$ lying in the linear space $(f^\perp)_s$: if there are any, then $f$ is the sum of $s$ codimension one forms, otherwise it is not. Using this we get a formula for $\smin$, thus solving the Waring Problem for codimension one decompositions.
\begin{thm}\label{sminTHM}
The generic form of degree $d$ in $n+1$ variables is the sum of $\min\{s: ns-{d-s+n\choose n }\geq 0\}$ codimension one forms and no fewer.
\end{thm}
\begin{proof}
The key part of the proof is Lemma \ref{Cod1ApolarityLemma}: a form $f\in S_d$ is the sum of $s$ codimension one forms if and only if $f^\perp$ contains a totally decomposable form of degree $s$ without repeated factors. Or, more geometrically, if and only if the linear space $F^\perp=(f^\perp)_s$ intersects $\Delta_{n,s}$ in at least a smooth point.
Consider the incidence correspondence
\[\Sigma=\{ (f,D) : D\in F^\perp, D=L_1\cdot\ldots\cdot L_s\}\subset\PP S_d\times \Delta_{n,s}\]
and the incidence diagram
\[\xymatrix{&{\Sigma}\ar[dl]_\alpha\ar[dr]^\beta&\\{\PP S_d} & &{\Delta_{n,s}} .}\]
Clearly $\dim\Sigma=\dim\PP S_d+ns-{d-s+n \choose n}$ (use $\beta$ to show that $\Sigma$ is a projective bundle over $\Delta_{n,s}$ having fiber over $D$ the projectivized of $(D^{-1})_d$). Moreover, for $s\geq\bar s=\min\{t: nt-{d-t+n \choose n}\geq 0\}$, the map $\alpha$ is surjective (a dimension count shows that $\Delta_{n,s}\cap F^\perp\neq\emptyset$ for {\em any} $f$). Let $\Sigma_0=\{(f,D=L_1\cdot\ldots\cdot L_s):D\in F^\perp, L_i\sim L_j \mbox{ for some  }i\neq j\}$ and notice that $\Sigma_0$ has codimension 2 in $\Sigma$. For $s\geq\bar s$ and $f$ generic, a dimension argument yields $\alpha^{-1}(f)\setminus\Sigma_0\neq\emptyset$, hence the claim.
\end{proof}
\begin{rem}\label{allPOSSdec}
Notice that, given $f\in S_d$, the variety $(f^\perp)_s\cap\Delta_{n,s}$ contains all the information about {\em all} the possible codimension one decompositions of the form $f$ involving $s$ summands. We will investigate this in the next section.
\end{rem}
\begin{rem}
In the case of binary forms ($n=1$), the Theorem gives $\smin(1,d)=\min\{s: 2s-d-1)\geq 0\}=\lceil{d+1\over 2}\rceil$ as was well known to Sylvester.
\end{rem}
As a folklore result, it is interesting to compare $\smin$ and $\sexp$, thus measuring how codimension one decompositions are defective. Compared to the sum of powers case the result is quite surprising.
\begin{cor}\label{AHsimeq}Let $n\geq 2$. The minimal number of summands appearing in the codimension one decompositions of the generic form of degree $d$ in $n+1$ variables is the expected one (i.e. $\smin(n,d)=\sexp(n,d)$) if and only if
\[d=2,3\mbox{ for any } n\geq 2\]
or
\[d=4,5,6\mbox{ and }8 \mbox{ for }n=2.\]
\end{cor}
\begin{proof} The proof is mainly an exercise in arithmetic. First we notice that the quantity $s^*(d)=\lceil{d+1\over 2}\rceil$ (see Remark \ref{sstar}) also satisfies the inequality $s^*(d)\geq\sexp(n,d)\;(\dagger)$. Then the result follows by studying the equalities $\smin(n,d)=s^*(d)\;(\star)$ and $\sexp(n,d)=s^*(d)\;(\sharp)$. The $n=2$ case is contained in \cite{Ca04}, thus we restrict our analysis to the $n>2$ cases.

As ${1\over n+{d+n-1\choose n-1}}{d+n\choose n}\leq{1\over n}(d+n)$, to show $(\dagger)$ it is enough to show that $s^*(d)\geq{1\over n}(d+n)$. For $d=2k$ this is equivalent to $(n-2)k\geq 0$. While, for $d=2k+1$, we get the inequality $(n-2)k\geq 1$. Thus $(\dagger)$ holds for $n> 2$ and $d>1$.

To study $(\star)$, we notice that $\smin(n,d)>s^*(d)$ if and only if $s^*(d)n-{d-s^*(d)+n \choose n}<0$. For $d=2k+1$ the last inequality is equivalent to $n(k+1)<{k+n\choose k}$, which holds for $k=2$ and can be proved to hold for $n>2$ and $k\geq 2$ by induction (use the fact that ${n+k\choose k+1}\geq n$). For $d=2k$ we have the inequality $n(k+1)<{k+n-1\choose k}$ which holds for $n>4$ and $k\geq 3$ and for $n=3,4$ and $k>3$. In conclusion, $(\star)$ could only possibly holds for $(n,d)=(3,6),(4,6),(n,2),(n,3),(n,4)$ for $n>2$ and it is easy to check that this is actually the case.

Finally, by direct substitution, we verify for which pairs $(\sharp)$ and $(\star)$ have common solutions.
\end{proof}

\section{How many decompositions?}\label{VSH}
In the previous section we solved the Waring Problem for codimension one decompositions, i.e. we determined the minimal number of summands appearing in the decompositions of a generic form. Once we know that a generic $f$ can be written as the sum of $\smin$ codimension one forms, it is natural to study in how many ways such a decomposition can be obtained. Actually, we are interested in an even more general question:
\begin{quote}
{\em How can we describe the codimension one decompositions of $f$ involving $s$ summands?}
\end{quote}
The answer is suggested by Remark \ref{allPOSSdec}:
\begin{defn}\label{defVS}
Let $F\subset\PP^n$ be a generic degree $d$ hypersurface, then the {\em variety of sums of codimension one forms} of $F$ with respect to $s$ is the scheme-theoretic intersection
\[\VS(F,s)=F^\perp\cap\Delta_{n,s}.\]
\end{defn}
\begin{rem}Recall that $F=V(f)$ for a generic $f\in S_d$ and that $F^\perp$ denotes the projectivization of the appropriate homogeneous piece of $f^\perp$. We define $\VS$ in terms of $F$ rather than of $f$ because the variety is a controvariant of the form under the action of $PGL(n+1)$.
\end{rem}
\begin{rem}Strictly speaking, $\VS(F=V(f),s)$ does not describe {\em all} the possible codimension one decompositions of $f$. Indeed, given a reduced point $L_1\cdot\ldots\cdot L_s\in \VS(F,s)$ we know that there exist codimension one forms $\hat f_1,\ldots,\hat f_s$ such that
 \[L_i\circ f_i=0, i=1\ldots s,\mbox{ and }f=\sum_1^s \hat f_i,\] but the forms $\hat f_i$'s are not uniquely determined by the $L_i$'s. It is not difficult to see that {\em all} the codimension one decompositions of $f$ are described by a projective bundle over $\VS(F,s)$.
\end{rem}

To carry on our analysis we need to study the variety $\Delta_{n,s}$ in some detail. Besides the general results in \cite{GKZ} (see Remark \ref{chowREM}), very few is known about $\Delta_{n,s}$: for $n=2$ and $d=3$ an invariant theory description of the defining ideal is contained in \cite{Chi02}; in \cite{Ca04JA} some geometric properties are described for $n=2$, any $d$, and we generalize the ideas contained there.
\begin{prop}\label{degDIMprop}
Let $F=L_1\cdot\ldots\cdot L_s\in\Delta_{n,s}$ be a generic point. If we let $X_F=\cup_{i\neq j}\{L_i=0\}\cap\{L_j=0\}\subset\check\PP^n$, then the tangent space to $\Delta_{n,s}$ in $F$ is
\[T_F(\Delta_{n,s})=|sH-X_F|=\{\mbox{hypersurfaces of degree $s$ in $\check\PP^n$ passing through $X_F$}\}.\]
Moreover
\[\begin{array}{rcl}
\deg\Delta_{n,s} & = & \#\{\mbox{$s$-uples of hyperplanes passing through $ns$ generic points in $\check\PP^n$}\}\\
                 & = & {ns-1\choose n-1}\cdot{n(s-1)-1\choose n-1}\cdot\ldots\cdot 1 .
\end{array}\]
\end{prop}
\begin{proof} Using the differential of a parametric description of $\Delta_{n,s}$ it is immediate to see that $T_F(\Delta_{n,s})$ is the projectivization of the vector space
\[\langle L_1\cdot\ldots \overline{L_i} \ldots\cdot L_s : \overline{L_i}\in T_1,i=1,\dots,s\rangle.\]
We claim that:
\begin{quote}
{\em the defining ideal of $X_F$, $I_{X_F}=\cap_{i\neq j}(L_i,L_j)$, is generated by the degree $s-1$ elements $\{L_1\cdot\ldots \widehat{L_i} \ldots\cdot L_s : i=1,\dots,s\}$.}
\end{quote}
This is enough to get the desired description of the tangent space.

To determine $\deg \Delta_{n,s}$ choose $ns$ generic points $P_1,\ldots,P_{ns}\in\check\PP^n$ and consider $H(P_1,\ldots,P_{ns})\subset\PP T_d$ which denotes the linear system of degree $s$ hypersurfaces through them. If $\YY=H(P_1,\ldots,P_{ns})\cap\Delta_{n,s}$ is zero dimensional and smooth, then its cardinality is the degree of the variety of totally decomposable forms. By the genericity of the $P_i$'s, $\YY$ is clearly a set of points of the desired cardinality. Notice that $G=R_1\cdot\ldots\cdot R_s\in\YY$ is a singular point if and only if
\[T_G(\Delta_{n,s})\cap H(P_1,\ldots,P_{ns})\neq G.\]
Or, in other terms, if and only if there are degree $s$ hypersurfaces passing through $X_G=\cup_{i\neq j}\{R_i=0\}\cap\{R_j=0\}$ and the $P_i$'s beside $\{G=0\}$. Notice that such an element meets the hyperplane $\{R_i=0\}\subset\check\PP^n$ in $s-1$ codimension 2 linear spaces and in $n$ points in generic position, thus $\{R_i=0\}$ is a component. Hence $G$ is the unique element with the required property and the degree formula is proved.
\begin{proof}[Proof of the claim]
The proof is purely algebraic. In the polynomial ring $\CC[Y_1,\ldots,Y_s]$, clearly we have
\[(Y_1\cdot\ldots\widehat{Y_i}\ldots\cdot Y_s: i=1,\ldots,s)=\cap_{i\neq j}(Y_i,Y_j)\]
and the ideal can be shown to be Cohen-Macaulay. By specializing we get
\[(L_1\cdot\ldots\widehat{L_i}\ldots\cdot L_s: i=1\ldots,s)\subseteq\cap_{i\neq j}(L_i,L_j)\]
where we know that the leftmost ideal is Cohen-Macaulay. Comparing degrees and dimensions the equality follows.
\end{proof}
\end{proof}
\begin{rem} This Proposition allows us to easily compute the degree of the varieties of secant lines to the quadratic Veronese varieties. In fact, it is well known that $\Sec_1(V_{n,2})=\Delta_{n,2}$ and hence $\deg \Sec_1(V_{n,2})={2n-1 \choose n-1}$.
\end{rem}

Using this result we can get some information on the varieties parameterizing codimension one decompositions:
\begin{thm}\label{DIMdegTHM}
Let $F\subset\PP^n$ be a generic degree $d$ hypersurface and let $s=\smin(n,d)$, then
\begin{itemize}
\item $\dim\VS(F,s)=ns-{d-s+n\choose n}$;
\item $\deg \VS(F,s)={ns-1\choose n-1}\cdot{n(s-1)-1\choose n-1}\cdot\ldots\cdot 1$.
\end{itemize}
\end{thm}
\begin{proof}With the notation of the Proof of Theorem \ref{sminTHM}, we have $\VS(F=V(f),s)=\alpha^{-1}(f)$. The dimension claim readily follows. Once $\dim\VS(F,s)$ is known, we realize that the intersection $F^\perp\cap\Delta_{n,s}$ is proper and hence we obtain the degree formula.
\end{proof}

The variety of totally decomposable forms contains singularities in codimension 2. Thus, as soon as $\dim\VS\geq 2$, the scheme parameterizing codimension one decomposition is singular and possibly not reduced. Nevertheless, we can get smoothness in one remarkable case:
\begin{prop}\label{PropRED}
Let $F\subset\PP^n$ be a generic degree $d$ hypersurface and let $s=\smin(n,d)$. If $\VS(F,s)$ is zero dimensional, then it is smooth (and hence reduced).
\end{prop}
\begin{proof}
The non-smoothness condition is a closed condition on the form $f\in S_d$ defining $F$. Hence we prove the Proposition by exhibiting forms with the required property. Notice that the zero dimension assumption yields the relation $ns={d-a+n\choose n}(\star)$.

Choose $ns$ generic linear forms $l_1,\ldots,l_{ns}$ and consider the form $g=\sum l_i^d$. Notice that $\VS(G=V(g),s)$ contains at least $\deg\Delta_{ns}$ points. Hence, {\em in the case} $\dim\VS(G,s)=0$, the variety is also smooth (the form $g$ is in general very degenerate as it is the sum of {\em few} powers, in particular Theorem \ref{DIMdegTHM} does not apply and the dimension has to be determined by other means). Let $\XX=\{l_1,\ldots,l_{ns}\}\subset\check\PP^n$ and denote by $I_{\XX}$ its defining ideal. Clearly $g^\perp\supset I_{\XX}$ by the Apolarity Lemma and, if equality holds in degree $s$, then $\dim\VS(G,s)=0$. By standard computation we get
\[\dim_\CC (I_{\XX})_s={s+n\choose n}-ns\mbox{ and } \dim_\CC (g^\perp)_s\geq{s+n\choose n}-{d-s+n\choose n}\]
and we want to show that the last inequality is an equality (this suffices by $(\star)$). If the inequality is strict, then there exists $D\in(g^\perp)_{d-s}$. Thus,
\[D\circ f=\sum_1^{ns}c_i D(l_i)l_i^s=0,\]
where $c_iD(l_i)\in\CC$ are not all zero as $(I_{\XX})_{d-s}=0$ by genericity. Hence a contradiction, as the $l_i$'s can be chosen in such a way that $l_1^d,\ldots,l_{ns}^s$ are linearly independent in $\PP S_s$ (notice that ${s+n\choose n}-ns=\dim_\CC (I_{\XX})_s>0$). In conclusion, $(g^\perp)_s= (I_{\XX})_s$ and $\dim\VS(G,s)=0$, as required.
\end{proof}
\begin{rem}
The previous Proposition has a ``negative'' consequence in the spirit of $19^{th}$ century invariant theory: there is no reasonable way (not even taking multiplicities in account) to produce canonical forms via codimension one decompositions. In other words, the decomposition is never unique (for similar results see,  for sum of powers, \cite{Me} and, for partially symmetric tensors, \cite{Fo}).
\end{rem}
\begin{rem} For $n=2$ it is easy to show that $\VS$ is zero dimensional in infinitely many cases ($\dim\VS=0$ is a degree two equation having integer solutions). If $n>2$, experiments suggests that zero dimensionality only occurs for $n=3,6$. Notice that the number of summands, $s$, rapidly increases and that the number of decompositions, $\#$, is soon intractable. This table shows what the numbers look like for $n\leq 100,d\leq 34$:
\[\begin{array}{llll}
d & n & s & \#
\\
\hline\\
5 & 2 & 3 & 15 \\
8 & 2 & 5 & 945 \\
17 & 6 & 14 & b\\
20 & 2 & 14 & 213458046676875 \\
25 & 2 & 18 & 221643095476699771875\\
34 & 3 & 28 & a\\
\end{array}\]
where $a$ and $b$ have 86 digits!
\end{rem}

\bibliographystyle{alpha}
\bibliography{MYbib}

\end{document}